\documentclass[11pt]{amsart}
\usepackage{amsmath,amssymb,mathrsfs}
\newtheorem{theorem}{Theorem}
\newtheorem{corollary}[theorem]{Corollary}
\newtheorem{lemma}[theorem]{Lemma}

\begin{document}

\title[Minimal surfaces in the unit ball]{A sharp bound for the area of minimal surfaces in the unit ball}
\author{Simon Brendle}
\address{Department of Mathematics \\ Stanford University \\ 450 Serra Mall, Bldg 380 \\ Stanford, CA 94305} 
\thanks{The author was supported in part by the National Science Foundation under grant DMS-0905628.}
\begin{abstract}
Let $\Sigma$ be a $k$-dimensional minimal surface in the unit ball $B^n$ which meets the boundary $\partial B^n$ orthogonally. We show that the area of $\Sigma$ is bounded from below by the volume of the unit ball in $\mathbb{R}^k$.
\end{abstract}
\maketitle

\section{Introduction}

One of the most important tools in minimal surface theory is the classical monotonicity formula, which asserts the following:

\begin{theorem}[cf. \cite{Allard}, \cite{Simon}]
\label{monotonicity.formula}
Let $\Sigma$ be a $k$-dimensional minimal submanifold of $\mathbb{R}^n$ with boundary $\partial \Sigma$. Moreover, let $y$ be a point in $\mathbb{R}^n$ and $r_0$ be a positive real number with the property that $\partial \Sigma \cap B_{r_0}(y) = \emptyset$. Then the function 
\[r \mapsto \frac{|\Sigma \cap B_r(y)|}{r^k}\] 
is monotone increasing for $r \in (0,r_0)$.
\end{theorem}

We note that Theorem \ref{monotonicity.formula} extends to the more general setting of stationary varifolds; see \cite{Allard}, Theorem 8.5. The monotonicity formula plays a fundamental role in the analysis of singularities. Moreover, it has a number of interesting geometric consequences (see e.g. \cite{Ekholm-White-Wienholtz}). In particular, it directly implies the following classical results: 

\begin{corollary}
\label{consequence.of.monotonicity.1}
Let $\Sigma$ be a $k$-dimensional minimal surface in the unit ball $B^n$ which passes through the origin and satisfies $\partial \Sigma \subset \partial B^n$. Then $|\Sigma \cap B^n| \geq |B^k|$.
\end{corollary}

Corollary \ref{consequence.of.monotonicity.1} can be viewed as a sharp version of F.~Almgren's density bound; see \cite{Almgren}, p.~343, for details.

\begin{corollary}
\label{consequence.of.monotonicity.2}
Let $\hat{\Sigma}$ be a closed minimal submanifold of the unit sphere $\partial B^n$ of dimension $k-1$. Then $|\hat{\Sigma}| \geq |\partial B^k|$.
\end{corollary}

In order to deduce Corollary \ref{consequence.of.monotonicity.2} from the monotonicity formula, one considers the $k$-dimensional minimal cone $\Sigma = \{\lambda \, x: x \in \hat{\Sigma}, \, \lambda>0\} \subset \mathbb{R}^n$. Even though the surface $\Sigma$ is singular at the origin, the monontonicity formula still holds, and we obtain 
\[|\hat{\Sigma}| = \lim_{r \to \infty} \frac{k \, |\Sigma \cap B_r(y)|}{r^k} \geq \lim_{r \to 0} \frac{k \, |\Sigma \cap B_r(y)|}{r^k} = k \, |B^k| = |\partial B^k|,\] 
where $y$ is arbitrary point on $\hat{\Sigma}$.

In \cite{Fraser-Schoen}, A.~Fraser and R.~Schoen considered a free boundary value problem for minimal surfaces in the unit ball. Specifically, they studied minimal surfaces in the unit ball which meet the boundary orthogonally. In this paper, we give an optimal lower bound for the area of such surfaces:

\begin{theorem}
\label{main.theorem}
Let $\Sigma$ be a $k$-dimensional minimal surface in the unit ball $B^n$. Moreover, suppose that the boundary of $\Sigma$ lies in the unit sphere $\partial B^n$ and meets $\partial B^n$ orthogonally. Then $|\Sigma| \geq |B^k|$. Moreover, if equality holds, then $\Sigma$ is contained in a $k$-dimensional subspace of $\mathbb{R}^n$.
\end{theorem} 

Applying the divergence theorem to the radial vector field $V(x) = x$ gives 
\[k \, |\Sigma| = \int_\Sigma \text{\rm div}_\Sigma V = \int_{\partial \Sigma} \langle V,x \rangle = |\partial \Sigma|.\] 
Hence, Theorem \ref{main.theorem} implies a sharp lower bound for the isoperimetric ratio of $\Sigma$: 

\begin{corollary}
Let $\Sigma$ be a $k$-dimensional minimal surface in the unit ball $B^n$. Moreover, suppose that the boundary of $\Sigma$ lies in the unit sphere $\partial B^n$ and meets $\partial B^n$ orthogonally. Then 
\[\frac{|\partial \Sigma|^k}{|\Sigma|^{k-1}} \geq \frac{|\partial B^k|^k}{|B^k|^{k-1}}.\] 
Moreover, if equality holds, then $\Sigma$ is contained in a $k$-dimensional subspace of $\mathbb{R}^n$.
\end{corollary} 

Theorem \ref{main.theorem} was conjectured by R.~Schoen (see e.g. \cite{Schoen}), following a question posed earlier by L.~Guth. The $k=2$ case of Theorem \ref{main.theorem} was verified by A.~Fraser and R.~Schoen (cf. \cite{Fraser-Schoen}, Theorem 5.4). 

The proof of Theorem \ref{main.theorem} is inspired by the classical monotonicity formula for minimal submanifolds, and its analogue for the mean curvature flow (cf. \cite{Ecker}, \cite{Huisken}). In order to prove the classical monotonicity formula, one applies the divergence theorem to the vector field $\frac{x-y}{|x-y|^k}$. This vector field can be interpreted as the gradient of the Newton potential in $\mathbb{R}^k$. Similarly, Huisken's monotonicity formula for the mean curvature flow involves an integral of the backward heat kernel in $\mathbb{R}^k$. In order to prove Theorem \ref{main.theorem}, we apply the divergence theorem to a suitably defined vector field $W$. This vector field agrees with the gradient of the Green's function for the Neumann boundary value problem on $B^k$, up to a factor.

The author would like to thank Professor Frank Morgan and Professor Brian White for comments on an earlier version of this paper.

\section{Proof of Theorem \ref{main.theorem}}

Let us fix a point $y \in \partial B^n$. We define a vector field $W$ in $B^n \setminus \{y\}$ by 
\[W(x) = \frac{1}{2} \, x - \frac{x-y}{|x-y|^k} - \frac{k-2}{2} \int_0^1 \frac{tx-y}{|tx-y|^k} \, dt.\] 

\begin{lemma}
\label{a}
For every point $x \in B^n$ and every orthonormal $k$-frame $\{e_1,\hdots,e_k\} \subset \mathbb{R}^n$, we have 
\[\sum_{i=1}^k \langle D_{e_i} W,e_i \rangle \leq \frac{k}{2}.\] 
\end{lemma}

\textbf{Proof.} 
We have 
\begin{align*} 
\sum_{i=1}^k \langle D_{e_i} W,e_i \rangle 
&= \frac{k}{2} - \frac{k}{|x-y|^{k+2}} \, \Big ( |x-y|^2 - \sum_{i=1}^k \langle x-y,e_i \rangle^2 \Big ) \\ 
&- \frac{k-2}{2} \int_0^1 \frac{tk}{|tx-y|^{k+2}} \, \Big ( |tx-y|^2 - \sum_{i=1}^k \langle tx-y,e_i \rangle^2 \Big ) \, dt \\ 
&\leq \frac{k}{2}, 
\end{align*} 
as claimed. \\

\begin{lemma}
\label{b}
The vector field $W$ is tangential along the boundary $\partial B^n$.
\end{lemma}

\textbf{Proof.} 
A straightforward computation gives 
\begin{align*} 
\langle W,x \rangle 
&= \frac{1}{2} \, |x|^2 - \frac{\langle x-y,x \rangle}{|x-y|^k} - \frac{k-2}{2} \int_0^1 \frac{\langle tx-y,x \rangle}{|tx-y|^k} \, dt \\ 
&= \frac{1}{2} \, |x|^2 - \frac{\langle x-y,x \rangle}{|x-y|^k} + \frac{1}{2} \int_0^1 \frac{d}{dt} \Big ( \frac{1}{|tx-y|^{k-2}} \Big ) \, dt \\ 
&= \frac{1}{2} \, |x|^2 - \frac{\langle x-y,x \rangle}{|x-y|^k} + \frac{1}{2} \, \frac{1}{|x-y|^{k-2}} - \frac{1}{2} \\ 
&= \frac{1}{2} \, (1-|x|^2) \, \Big ( \frac{1}{|x-y|^k} - 1 \Big ). 
\end{align*} 
In particular, $\langle W,x \rangle = 0$ for all $x \in \partial B^n \setminus \{y\}$. \\

\begin{lemma}
\label{c}
We have 
\[W(x) = -\frac{x-y}{|x-y|^k} + o \Big ( \frac{1}{|x-y|^{k-1}} \Big )\] 
as $x \to y$.
\end{lemma} 

\textbf{Proof.} 
Using the inequality 
\[|tx-y|^2 = t \, |x-y|^2 + (1-t)(1-t \, |x|^2) \geq t \, |x-y|^2 + (1-t)^2,\] 
we obtain 
\begin{align*} 
\bigg | W(x) + \frac{x-y}{|x-y|^k} \bigg | 
&= \bigg | \frac{1}{2} \, x - \frac{k-2}{2} \int_0^1 \frac{tx-y}{|tx-y|^k} \, dt \bigg | \\ 
&\leq \frac{1}{2} + \frac{k-2}{2} \int_0^1 \frac{1}{|tx-y|^{k-1}} \, dt \\ 
&\leq \frac{1}{2} + \frac{k-2}{2} \int_0^1 \Big ( \frac{1}{t \, |x-y|^2 + (1-t)^2} \Big )^{\frac{k-1}{2}} \, dt. 
\end{align*} 
It follows from the dominated convergence theorem that 
\[\int_0^1 \Big ( \frac{|x-y|^2}{t \, |x-y|^2 + (1-t)^2} \Big )^{\frac{k-1}{2}} \, dt \to 0\] 
as $x \to y$. Thus, we conclude that 
\[W(x) + \frac{x-y}{|x-y|^k} = o \Big ( \frac{1}{|x-y|^{k-1}} \Big )\] 
as $x \to y$. This completes the proof. \\

We now give the proof of Theorem \ref{main.theorem}. To that end, let us fix a point $y \in \partial \Sigma$, and let $W$ be the vector field defined above. Note that $W$ is smooth on $B^n \setminus \{y\}$. Using the divergence theorem, we obtain 
\begin{align} 
\label{first.variation}
&\int_{\Sigma \setminus B_r(y)} \Big ( \frac{k}{2} - \text{\rm div}_\Sigma W \Big ) \notag \\ 
&= \frac{k}{2} \, |\Sigma \setminus B_r(y)| - \int_{\Sigma \cap \partial B_r(y)} \langle W,\nu \rangle - \int_{\partial \Sigma \setminus B_r(y)} \langle W,x \rangle, 
\end{align}
where $\nu$ denotes the inward pointing unit normal to the region $\Sigma \cap B_r(y)$ within the submanifold $\Sigma$. In particular, the vector $\nu$ is tangential to $\Sigma$, but normal to $\Sigma \cap \partial B_r(y)$. It is easy to see that 
\[\nu = -\frac{x-y}{|x-y|} + o(1)\] 
for $x \in \Sigma \cap \partial B_r(y)$. Using Lemma \ref{c}, we obtain 
\[\langle W,\nu \rangle = \frac{1}{r^{k-1}} + o \Big  ( \frac{1}{r^{k-1}} \Big )\] 
for $x \in \Sigma \cap \partial B_r(y)$. Since 
\[|\Sigma \cap \partial B_r(y)| = \frac{1}{2} \, |\partial B^k| \, r^{k-1} + o(r^{k-1}),\] 
we conclude that 
\begin{equation} 
\label{term.1}
\lim_{r \to 0} \int_{\Sigma \cap \partial B_r(y)} \langle W,\nu \rangle = \frac{1}{2} \, |\partial B^k| = \frac{k}{2} \, |B^k|. 
\end{equation}
On the other hand, it follows from Lemma \ref{b} that $\langle W,x \rangle = 0$ for all $x \in \partial B^n \setminus \{y\}$. This implies 
\begin{equation} 
\label{term.2}
\int_{\partial \Sigma \setminus B_r(y)} \langle W,x \rangle = 0. 
\end{equation}
Combining (\ref{first.variation}), (\ref{term.1}), and (\ref{term.2}), we obtain 
\[\lim_{r \to 0} \int_{\Sigma \setminus B_r(y)} \Big ( \frac{k}{2} - \text{\rm div}_\Sigma W \Big ) = \frac{k}{2} \, (|\Sigma| - |B^k|).\] 
On the other hand, we have 
\[\frac{k}{2} - \text{\rm div}_\Sigma W \geq 0\] 
by Lemma \ref{a}. Putting these facts together, we conclude that $|\Sigma| - |B^k| \geq 0$. 

It remains to analyze the case of equality. Suppose that $|\Sigma| - |B^k| = 0$. This implies 
\[\frac{k}{2} - \text{\rm div}_\Sigma W = 0\] 
at each point $x \in \Sigma \setminus \{y\}$. Consequently, for each point $x \in \Sigma \setminus \{y\}$, we have 
\[|x-y|^2 - \sum_{i=1}^k \langle x-y,e_i \rangle^2 = 0,\] 
where $\{e_1,\hdots,e_k\}$ is an orthonormal basis of $T_x \Sigma$. From this, we deduce that $x-y \in T_x \Sigma$ for all points $x \in \Sigma \setminus \{y\}$. Since $y \in \partial \Sigma$ is arbitrary, we conclude that $\partial \Sigma$ is contained in a $k$-dimensional affine subspace of $\mathbb{R}^n$. By the maximum principle, $\Sigma$ is contained in a $k$-dimensional affine subspace of $\mathbb{R}^n$.

\end{document}